\documentclass[10pt]{article}
\usepackage{amsfonts}
\usepackage{amssymb}
\usepackage{amsthm}
\usepackage{amsmath} 
\usepackage{graphicx}
\usepackage{color}
\usepackage[mathscr]{euscript}
\date{}
\newcommand{\sss}{\setcounter{equation}{0}}
\newtheorem{theorem}{THEOREM}[section]

\newtheorem{lemma}[theorem]{LEMMA}

\newtheorem{remark}[theorem]{REMARK}

\newtheorem{definition}[theorem]{DEFINITION}

\newcommand\beq{\begin{equation}}
\newcommand\ene{\end{equation}}
\def \ds{\displaystyle}

\begin{document}
\baselineskip=12 pt
\parskip 6 pt

\title{ Trace maps   under weak regularity assumptions. \thanks{  Research partially supported by project  PAPIIT-DGAPA UNAM  IN103918 and by project SEP-CONACYT CB 2015, 254062. }}
\author{ Ricardo Weder\thanks {Fellow, Sistema Nacional de Investigadores.}\thanks{ Electronic mail: weder@unam.mx. Home page: http://www.iimas.unam.mx/rweder/rweder.html} \\
Departamento de F\'{\i}sica Matem\'atica.\\
 Instituto de Investigaciones en Matem\'aticas Aplicadas y en
 Sistemas. \\
 Universidad Nacional Aut\'onoma de M\'exico.\\
  Apartado Postal 20-126,
Ciudad de M\'exico  01000, M\'exico.}

\maketitle


\centerline{\it In Memory of Erik Balslev}

\vspace{.5cm}
 \centerline{{\bf Abstract}}
 We study  bounded trace maps  on {\it hypersurfaces } for Sobolev spaces   from a point of view that is fundamentally different from the one in the classical theory. This allows us to construct   bounded trace maps under weak regularity assumptions on the {\it hypersurfaces}. In the case of bounded domains  in $\mathbf R^n$ we only require the continuity of the boundary. For  {\it hypersurfaces} in the whole space $\mathbf R^n$ we only assume that the  {\it hypersurfaces} are Lebesgue measurable.
 As an application of our trace maps we consider the Dirichlet  problem and we prove a coarea formula where the level sets are only assumed to be Lebesgue measurable  {\it hypersurfaces}.


\bigskip

\section{Introduction} \sss
The study of trace maps on   {\it hypersurfaces},  for Sobolev spaces, is a fundamental problem  that has been considered extensively. It plays an essential role in many areas of analysis. For example,  it is crucial in the theory of Sobolev spaces and  in the formulation and  the solution of boundary value problems. See for example  \cite{24}, \cite{ah},  \cite{ho}, \cite{jk},\cite{lm}, \cite{necas}, and  \cite{rs2}. The more general results in these  references  require that the {\it  hypersurface} is $C^1$ or Lipschitz continuous.   This is a strong  restriction. Note that the derivatives of Lipschitz continuous functions are bounded and that this is a condition  that often is not satisfied in the applications. Moreover,  in \cite{mp} there are  results on trace maps in non-Lipschitz domains bounded by Lipschitz surfaces, in domains with cusps and with peaks and in capacity criteria for trace maps.   Furthermore, in \cite{m} there are results on trace maps for function of bounded variation on domains with finite perimeter and with a normal vector at the boundary. Moreover,   \cite{m} gives capacity criteria for the existence of trace maps. Furthermore, \cite{bk} gives  trace maps for functions of bounded variation, assuming that the boundary  is a $(n-1)$-rectifiable set and that it has a normal vector. For further results see \cite{jw} and \cite{sh}. 

We consider {\it hypersurfaces} that are star shaped about the origin, but we impose weak  conditions on the smoothness of the {\it hypersurfaces}. In the case when the {\it hypersurface} is the boundary of a bounded star shaped domain we only require that the function that characterizes the {\it hypersurface}  is continuos. Further, when  we take the trace of functions in Sobolev spaces in $\mathbf R^n$  in a star shaped {\it hypersurface} we only require that the function that characterizes the {\it hypersurface} is Lebesgue measurable. Star shaped are a special type of  {\it hypersurfaces}, but note that this condition is often satisfied in the applications.

The paper is organized as follows. In Section 2 we construct trace maps in the boundary of bounded star shaped domains, and we apply our results to  the Dirichlet problem. In Section 3 we construct  trace maps in star shaped {\it hypersurfaces}, and we discuss the relation   of our results  with our previous results in \cite{we}, where we constructed trace maps in the slowness surface of strongly propagative systems of  equations. The trace maps in \cite{we} were applied to the spectral and the scattering theory of strongly propagative systems of equations. In Section 4 we consider the relation  of our trace maps of Section 3 with the coarea formula.
 We prove a coarea formula where the level sets are star shaped  {\it   hypersurfaces} that are  only assumed to be Lebesgue measurable.
 Note that  using the the usual localization arguments, with partitions of unity, we can extend our results to domains that are locally star shaped.

\section{Domains with  star shaped boundary}\label{sec2}
\sss

We first introduce some standard definitions and notations. For any open set $\Omega \subset \mathbf R^n$ we denote by $C^1(\Omega)$ the set of all continuously differentiable functions in $\Omega$, and by $C^1(\overline{\Omega})$ the set of all functions in $C^1(\Omega)$ that together with all its first derivatives have continuous extensions to $\overline{\Omega}.$ Furthermore, we designate by 
 $C^\infty(\Omega)$ the set of all infinitely 
differentiable functions in $\Omega,$  by $C^\infty_0\Omega)$ the set of all infinitely differentiable functions with compact support in $\Omega$ and by $C^\infty(\overline{\Omega})$ the set of all functions in $C^\infty(\Omega)$ that together with all its derivatives have a continuous extension to $\overline{\Omega}$. By $C^\infty_0(\mathbf R^n)\big|_{\Omega}$   we denote the set of all the restrictions to $\Omega$ of functions in $C^\infty_0(\mathbf R^n).$ We denote by $\mathcal S$ the space of Schwartz of all infinitely differentiable function on $\mathbf R^n$ that  together with all its derivatives remain bounded when multiplied by any polynomial. We  denote by $C$ a generic constant that can take different values when it appears in various places.
The symbol  $L^p(\Omega), 1 \leq p < \infty$ denotes the  standard   space of Lebesgue measurable functions in $\Omega$ whose absolute value to the power $p$ is integrable, with norm,
$$
\|f\|_{L^p(\Omega)}:= \left( \int_\Omega |f(x)|^p \, dx \right)^{1/p}.
$$
$L^2(\Omega)$ is a Hilbert space with the standard scalar product,
$$
\left( f,g \right)_{L^2(\Omega)}:= \left(  \int_{\Omega}\, f(x)\, \overline{g(x)}\, dx      \right)^{1/2}.
$$
 
The Sobolev space $W^{(1)}_p(\Omega), 1 \leq p  < \infty,$ is the Banach  space of all functions $f(x) \in L^p(\Omega)$ such that the derivatives in distribution sense $ \frac{\partial}{\partial x_i} f(x), 1 \leq i \leq n,$ are functions in $L^p(\Omega).$
The norm of $W^{(1)}_p(\Omega)$ is given by,
\beq\label{2.3}
\|f\|_{W^{(1)}_p(\Omega)}:= \left( \int_\Omega |f(x)|^p \, dx+\sum_{i=1}^n \int_\Omega \left|\frac{\partial}{\partial x_i}f(x)\right|^p \, dx  \right)^{1/p}.
\ene

We denote by $ W^{(1)}_{p,0}(\Omega)$ the closure of $C^\infty_0(\Omega)$ in the norm of $W^{(1)}_p(\Omega).$

A domain  $\Omega$ is an  open set in $ \mathbf R^n, n \geq 2.$ Let $ \partial \Omega$ be its boundary.   Let us denote by $S_{n-1}$ the unit sphere in $\mathbf R^n$.  We consider domains that are star shaped with respect to the origin. They are defined as follows.
\begin{definition} \label{def2.1}
{\rm The  domain $\Omega$ is star shaped with respect to the origin with a continuos function that characterizes the boundary, if there is a continuous function $b(\nu) >0 $ defined for $ \nu \in S_{n-1}$ such that,}
\beq\label{2.1} 
\Omega=\left\{ x \in \mathbf R^n: |x| < b\left(\frac{x}{|x|}\right)\right\}.
\ene
\end{definition}

Note that
\beq\label{2.2}
\partial \Omega=\left\{x \in \mathbf R^n: |x|= b\left(\frac{x}{|x|}\right)  \right\}.
\ene
 In this section we always assume that $\Omega$ is star shaped with respect to the origin with a continuous function that characterizes the boundary.  By Theorem 3.2 in page 67 of \cite{necas}, or by Theorem 1 in page 10 of \cite{m}, if $ \Omega$ is  star shaped  with respect to the origin with a continuous boundary, $C^\infty_0(\mathbf R^n)\big|_{\Omega}$ is dense in $W^{(1)}_p(\Omega)$ and as  $C_0^\infty(\mathbf R^n)\big|_{\Omega} \subset C^\infty(\overline{\Omega}),$ also $C^\infty(\overline{\Omega}) $ is dense in  $W^{(1)}_p(\Omega).$

 In the classical theory the  trace maps are defined as bounded operators from $W^{(1)}_p(\Omega), 1 \leq p < \infty,$ into $L^p(\partial \Omega),$ where $L^p(\partial \Omega)$ is the standard Banach space of all Lebesgue-measurable functions  in $\partial\Omega$ whose absolute value to the power $p$ is Lebesgue integrable on $\partial \Omega.$ This requires that  $\partial \Omega$ is Lipschitz. To understand the origin of this restriction it is instructive to consider the case of domains in two dimensions.  Let us use polar coordinates in $\mathbf R^2,$ $x_1= r\, \cos\theta, x_2=  r \, \sin\theta,$ with $r:= |x|, 0 \leq \theta < 2\pi.$ In these coordinates the unit vector $\nu$ in $S_{n-1}$ can be written as follows
 \beq\label{unitvec}
 \nu= \nu(\theta)= \cos\theta \,e_{x_1} + \sin\theta\, e_{x_2},
 \ene
 with $e_{x_1},$  $e_{x_2},$ respectively, the unit vectors along the $x_1$ and the $x_2$ axis. 
Moreover, for simplicity we denote,
$$
b(\theta):= b(\nu(\theta)).
$$
Furthermore, by \eqref{2.2}  $\partial\Omega$  is given in parametric form in the following way,
$$
\partial \Omega= \left\{  x \in \mathbf R^2: x =\varphi(\theta), 0 \leq \theta < 2 \pi \right\},
$$
where
$$
\varphi(\theta):= b(\theta)\, \nu(\theta),\quad 0 \leq \theta < 2 \pi.
$$
Assume  that $\partial \Omega$ is Lipschitz. For any  set $O \subset \partial \Omega,$ that is Lebesgue measurable, let us  denote by $m_{L,S}(O)$  its  Lebesgue measure.
Then, for any function  $f \in L^p(\partial \Omega)$ the integral over $\partial \Omega$ of its absolute  value to the power $p$ is given by, 
\beq\label{2.2.b}
\int_S\, |f(\omega)|^p\, dm_{L,S}(\omega)= \int_0^{2\pi}\,  |f( b(\theta)\, \nu(\theta))|^p\,  |\varphi'(\theta)|\,  d\theta.
\ene

Equation \eqref{2.2.b} shows clearly why in the classical theory of trace maps from $W^{(1)}_p(\Omega)$ into $L^p(\partial \Omega)$ it is necessary to impose the Lipschitz condition on $\partial \Omega.$ Namely, if  $   \varphi(\theta)$ is Lipschitz, we have that $ \varphi'(\theta)$ exists for almost every $\theta$ and it is bounded. In consequence, the right-hand side of \eqref{2.2.b} is well defined. In our case we only assume that $\varphi(\theta)$ is a continuous function,  that does not have to be differentiable. Hence, in our situation the right-hand side of \eqref{2.2.b}  makes no sense. So, in our case there no hope of constructing bounded trace maps with target space the standard $L^p(\partial\Omega),$ and we have to look for a different target space. This is certainly a radical departure from the standard point of view. There is, of course, a long tradition of taking the standard $L^p(\partial\Omega)$ as target space. However, on spite of this, there is no fundamental reason to take as target space the standard $L^p(\partial\Omega).$   What is important is to obtain trace maps that are useful tools in applications, for example, to boundary value problems,  to spectral and scattering theory, or to a coarea formula.  So, what are the options that we have in our case where $\varphi(\theta)$ is only assumed to be continuous. The simplest solution to our problem is just to eliminate  $ |\varphi'(\theta)|$  from the right-hand side of \eqref{2.2.b} and to propose for the integral of $ |f(\omega)|^p$ over $S$ the quantity,
$$
 \int_0^{2\pi}\,  |f( b(\theta)\, \nu(\theta))|^p\,  d\theta.
$$
Going back to the $n$ dimensional case, for a function $f(\omega), \omega \in S$ we propose for the integral of $|f(\omega)|^p$ over $S$ the quantity
\beq\label{2.2.c}
\int_{S_{n-1}}\,|f(b(\nu)\,\nu |^p\, dm_{S_{n-1}}(\nu),
\ene
where  for any set $A$ that is Lebesgue measurable in $S_{n-1},$ by $m_{S_{n-1}}(A)$ we denote its Lebesgue measure. It can be said that in \eqref{2.2.c} we take the pull back of the Lebesgue measure under the projection of $\partial \Omega$ onto $S_{n-1}.$ We show below that with our proposal we obtain trace maps that have the main properties that make the standard trace maps useful, when the later exist.  
To accomplish this task we first have to precisely define our target spaces for the trace maps. 
    
We proceed to define a measure in $\partial \Omega$ that is appropriate for our purposes. We define the following bijective  function, $\mathcal P,$ from $\partial \Omega$ onto $S_{n-1}.$ 
\beq\label{2.4}
\mathcal P(\omega):= \nu= \frac{\omega}{|\omega|}, \,\text{for}\, \omega \in \partial \Omega.
\ene 
The inverse function is given by,
\beq\label{2.5}
\mathcal P^{-1}(\nu):= \omega= b(\nu)\, \nu, \, \text{for}\, \nu \in S_{n-1}.
\ene
Clearly $\mathcal P$ is onto and $ \mathcal P$ and $\mathcal  P^{-1}$ are continuous. As usual for $ O \subset \partial \Omega,$ we denote, $ \mathcal P(O):=\{ \nu \in S_{n-1}: \nu= \mathcal P(\omega), \, \text{for some}\, \omega \in O\}.$ 
\begin{definition}\label{def2.2}{\rm
We say that a set $ O \subset \partial \Omega$ is $\partial \Omega$-measurable if $ \mathcal P(O)$ is Lebesgue measurable in $S_{n-1}.$ Further, for any  
$\partial \Omega$-measurable set $O,$ we define its measure, in symbols $ m_{\partial \Omega}(O),$ by
\beq\label{2.6} 
m_{\partial \Omega}(O):= m_{S_{n-1}}(\mathcal P(O)).
\ene
 We designate by $\mathcal M_{\partial \Omega}$ the set of all $\partial \Omega$-measurable sets. Clearly,  $\mathcal M_{\partial \Omega}$ is a $\sigma$-algebra and $m_{\partial \Omega}$ is a $\sigma$-additive measure on   $\mathcal M_{\partial \Omega}.$}
\end{definition}
Let $f(\omega)$ be a function defined on $\partial\Omega.$ We denote by $f_{S_{n-1}}(\nu)$ the function,
\beq\label{2.7}
f_{S_{n-1}}(\nu):= f(b(\nu)\,  \nu), \, \nu \in S_{n-1}.
\ene
Note that for any set $ A \subset \mathbf R,$
$$
\mathcal P\left(f^{-1}(A)\right)= f^{-1}_{S_{n-1}}( A).
$$
Observe  that  $f(\omega)$ is $\partial \Omega$-measurable if and only if $f_{S_{n-1}}(\nu)$ is Lebesgue measurable.
For $ \partial \Omega$-integrable functions $f(\omega),$ we have that,
\beq\label{2.9}
\int_{\partial \Omega}\, f(\omega)\, dm_{\partial \Omega}(\omega)= \int_{S_{n-1}}\,f_{S_{n-1}}(\nu)\, dm_{S_{n-1}}(\nu).
\ene

\begin{definition} \label{def2.3}\rm
We denote by $\mathcal L^p(\partial \Omega), 1 \leq p < \infty,$ the Banach space of all complex valued, $\partial \Omega$-measurable functions, $ f(\omega),$ such that
$|f(\omega)|^p$  is  integrable over $\partial \Omega$ with respect to the measure $m_{\partial \Omega}.$ The norm  in   $\mathcal L^p(\partial \Omega)$  is given by
\beq\label{2.10}
\left\| f \right\|_{\mathcal L^p(\partial \Omega)}:=\left( \int_{\partial \Omega} |f(\omega)|^p\, dm_{\partial \Omega}(\omega)   \right)^{1/p}.
\ene
\end{definition}
We have that,
\beq\label{2.11}
\left\|\ f \right\|_{\mathcal L^p(\partial \Omega)}=\left( \int_{S_{n-1}}\, |f_{S_{n-1}}(\nu)|^p\, dm_{S_{n-1}}(\nu)   \right)^{1/p}.
\ene

 
 
  We now state our trace map theorem.
\begin{theorem}\label{theo2.4}
Suppose that the domain $\Omega$ is star shaped with respect to the origin with a continuos function that characterizes the boundary. Then, there is a trace map $T_p$ that is bounded from $W^{(1)}_p(\Omega)$ into $\mathcal L^p(\partial \Omega), 1 \leq p < \infty,$ such that
\beq\label{2.12}
(T_pf)(\omega)= f(\omega), \, f(\omega)\in C^1(\overline{\Omega} ).
\ene
Furthermore, the range of $T_p$ is dense in $ \mathcal L^p(\partial \Omega), $
\beq\label{2.12.b}
\overline{T_p W^{(1)}_p(\Omega)}= \mathcal L^p(\Omega), 1 \leq p < \infty.
\ene
\end{theorem}

\noindent{\it Proof:} Assume that $ f(x)\in  C^1(\overline{\Omega} ).$  Let $ C_1 >0 $ be such that $ 0 < C_1 < b(\nu) < 1/ C_1, \nu \in S_{n-1}.$ Take $ h(r) \in C^\infty_0((0,\infty)), $ such that $ h(r)=0, 0< r <  C_1 /4, h(r)=1,  C_1/2 < r < 2/ C_1, h(r)=0, 
r >  3/ C_{1}.$ Denote $g(x)=  h(|x|)\,f(x).$ Then,
\beq\label{2.13}
f(b(\nu) \,\nu)=  g(b(\nu)\nu)=\int_0^{b(\nu)}\, \frac{\partial}{\partial \mu } g(\mu\, \nu)\, d\mu.
\ene
Hence, by the H\"older inequality,
 \beq\label{2.14}\begin{array}{l}
|f(b(\nu)\nu)|^p \leq C   \int_0^{b(\nu)}\,  \left|\frac{\partial}{\partial \mu } g(\mu\, \nu)\right|^p\, d\mu\leq \\\\
 C   \int_0^{b(\nu)}\,  \left|\frac{\partial}{\partial \mu } g(\mu\, \nu)\right|^p\,  \mu^{n-1}d\mu,
\end{array} \ene
 where in the last inequality we used that  as $ g(x)=0,$ for  $|x| \leq C_1/4,$ we have that $\frac{1}{\mu}\leq C$ on the support of $g(x).$ Moreover, by \eqref{2.9} and \eqref{2.14}
 \beq\label{2.15}\begin{array}{l}
\int_{\partial \Omega}\, |f(\omega)|^p\, dm_{\partial \Omega}(\omega)=
\int_{S_{n-1}}\, | f(b(\nu)\, \nu) |^p\,  dm_{S_{n-1}}(\nu)\, \leq \\\\
C \,\int_{S_{n-1}}\,  dm_{S_{n-1}}(\nu)\, \int_0^{b(\nu)}\,  \left|\frac{\partial}{\partial \mu } g(\mu\, \nu)\right|^p\,  \mu^{n-1}d\mu \leq \\\\
C \,\int_{S_{n-1}}\, dm_{S_{n-1}}(\nu)\, \int_0^{b(\nu)}\, \left|(\nabla g)(\mu\,\nu)\right|^p\,  \mu^{n-1}d\mu= \\\\
C\, \int_{\Omega}\, | \nabla g(x)|^p\, dx.
 \end{array}
\ene
But since $g(x)= h(|x|)\, f(x),$ by \eqref{2.15},
\beq\label{2.16}
\|f\|_{\mathcal L^p(\partial \Omega)} \leq C\, \|f\|_{W^{(1)}_p(\Omega)}.
\ene
Finally, the existence of $T_p$ follows since as $C^\infty(\overline{\Omega}) \subset C^1(\overline{\Omega}),$  we have that $C^1(\overline{\Omega})$ is dense in $W^{(1)}_p(\Omega).$

Let us now prove \eqref{2.12.b}. Suppose that $f(\omega) \in\mathcal  L^p(\partial \Omega).$ Then, $f_{S_{n-1}}(\nu):= f(b(\nu) \nu) \in L^p(S_{n-1}).$ Since $C^\infty(S_{n-1})$ is dense on $L^p(S_{n-1})$ there is sequence $g_m(\nu) \in C^\infty(S_{n-1})$ such that,
\beq\label{l2.a.1}
\lim_{m\to \infty}\| f_{S_{n-1}}-g_m \|_{L^p(S_{n-1})}=0.
\ene
Let us designate,
\beq\label{2.a.2}
f_m(x):= h(|x|)\, g_m\left(\frac{x}{|x|}\right) \in C^\infty_0(\mathbf R^n).
\ene
Denote by $f_m\big|_{\Omega}(x)$ the restriction of $f_m(x)$ to $\Omega.$ Then,  $f_m\big|_{\Omega}(x) \in C^\infty(\overline{\Omega}).$ Further,
\beq\label{2.a.3}  \begin{array}{l}
\lim_{m \to \infty}\|f -T_p  f_m\big|_{\Omega}  \|^p_{\mathcal L^p(\partial \Omega)}=\lim_{m \to \infty} \int_{S_{n-1}}\, |f(b(\nu)\nu)- h(b(\nu))\,  g_m(\nu)|^p\, dm_{S_{n-1}}(\nu)= \\\\
 \lim_{m \to \infty} \int_{S_{n-1}}\, |f(b(\nu)\nu)- g_m(\nu)|^p\, dm_{S_{n-1}}(\nu )=\lim_{m \to \infty} \int_{S_{n-1}}\, |f_{S_{n-1}}(\nu)- g_m(\nu)|^p\, dm_{S_{n-1}}(\nu)=0.
\end{array} \ene
\qed

Note that the proof of Theorem \ref{theo2.4} only requires that $\frac{\partial}{\partial \mu} f(\mu\nu) \in L^p(\Omega\setminus B_{C/4}),$ where $B_{C/4}$ is the ball of center zero and radius $C/4.$ This allows for a more general formulation of Theorem \ref{theo2.4}. 

The next theorem shows that our trace map $T_p$ characterizes $W^{(1)}_{p,0}(\Omega)$  as the Banach space of all functions in $W^{(1)}_p(\Omega)$ that are zero on $\partial\Omega$ in trace sense.
\begin{theorem}\label{theo2.5}
Suppose that the domain $\Omega$ is star shaped with respect to the origin with a continuos function that characterizes the boundary. Then, for $1 \leq p < \infty,$
\beq\label{2.16.b}
W^{(1)}_{p,0}(\Omega)= \left\{ f(\omega) \in W^{(1)}_p(\Omega): (T_pf)(\omega)=0 \right\}.
\ene
\end{theorem}
\noindent{\it Proof:} Clearly all the functions in $W^{(1)}_{p,0}(\Omega)$ are zero on $\partial \Omega$ in trace sense. Suppose that $ f(x)\in W^{(1)}_p(\Omega)$ and that 
$(T_p f)(\omega)=0.$  Then, eventually after redefining $f(x)$ as equal to zero in a set of  measure zero, we have that $f(\omega)=0, \omega \in \partial \Omega.$ 
By  Theorem 2.2 in page 61 and Theorem 2.3 in page 63 of \cite{necas}, or Theorem 1 in page four of \cite{m}, (see also \cite{le}), after redefining, if necessary, $f(x)$ as equal to zero in a set of measure zero, we can assume that  $f(x)$ is absolutely continuous on  all the lines parallel to the $x_1$ axis,  and that the classical derivative,  $[\frac{\partial}{\partial x_1}f(x)],$ coincides with the distribution derivative, i.e., $[\frac{\partial}{\partial x_1}f(x)]= \frac{\partial}{\partial x_1}f(x) \in L^p(\Omega).$ 
We extend $f(x)$ by zero to a function defined in $\mathbf R^n$ by setting $ f(x)=0, x \in \mathbf R^n\setminus \Omega$ and we also denote by $f(x)$ the  function extended to $\mathbf R^n.$
We proceed to prove that $\frac{\partial}{\partial x_1} f(x) \in L^p(\mathbf R^n).$  For any $(x_2,\cdots, x_n)\in \mathbf R^{n-1}$ we denote by 
$ L(x_2,\cdots,x_n)$ the line parallel to the $x_1$ axis that passes through the point $(0,x_2,\cdots,x_n),$
i.e.,
$$
L(x_2,\cdots,x_n)=\{ (x_1,x_2,\cdots, x_n), x_1 \in \mathbf R\}.
$$
Take any $\varphi(x)\in C^\infty_0(\mathbf R^n).$ Since $f(x)$ is absolutely continuous on the line $L(x_2, \cdots,x_n)$ and $f(x)=0,$ for $x \in \partial \Omega \cap L(x_2,\cdots, x_n),$ integrating by parts, we obtain that,   
\beq\label{2.16.c}\begin{array}{l}
\ds\int_{\mathbf R^n}\, f(x)\, \frac{\partial}{\partial x_1 }\varphi(x)\, d^nx = \int_{\mathbf R^{n-1}}\, dx_2 \cdots dx_n \int_{L(x_2,\cdots, x_n)}\, dx_1 f(x_1, x_2,\cdots, x_n)\, \frac{\partial}{\partial x_1}\varphi (x_1, x_2,\cdots,x_n)= \\\\
\ds-   \int_{\mathbf R^{n-1}} dx_2 \cdots dx_n \int_{L(x_2,\cdots, x_n)}\, dx_1  \left[\frac{\partial}{\partial x_1 }f(x_1, x_2,\cdots, x_n)\right]\, \varphi (x_1, x_2,\cdots,x_n)
\ds  =-\int_{\mathbf R^n}\,  \left[\frac{\partial}{\partial x_1 } f(x)\right]  \varphi(x)\, dx. 
\end{array}
\ene
Then, $\frac{\partial}{\partial x_1}f(x)= [\frac{\partial}{\partial x_1}f(x)]\in L^p(\mathbf R^n).$ We prove in the same way that $\frac{\partial}{\partial x_i}f(x)\in L^p(\mathbf R^n), i=2,\cdots,n.$ Hence $ f(x)\in W^{(1)}_p(\mathbf R^n).$ For $   \lambda > 1 $ we define  $ f_\lambda(x):=f(\lambda x).$ Then, $f_\lambda(x)\in W^{(1)}_p(\mathbf R^n)$ and 
$$
\ds \lim_{\lambda \to 1} \|f(x)-f_\lambda(x) \|_{W^{(1)}_{p}(\mathbf R^n)}=0.
$$
Finally, as for any fixed $\lambda >1,$ we have that $f_\lambda(x) =0,$ for $x \in \Omega \setminus \Omega_\lambda,$ where $\Omega_\lambda:= \{  x= r \nu, \nu \in S_{n-1},  0 \leq  r <\frac{b(\nu)}{\lambda} \},$  and $\overline{\Omega_\lambda} \subset \Omega$ it follows that $f_\lambda(x) \in W^{(1)}_{p,0}(\Omega).$
\qed

By theorems \ref{theo2.4} and \ref{theo2.5} our trace map $T_p$ have the following properties,
\begin{enumerate}
\item
For every $ f \in C^1(\overline{\Omega}),\, T_p(f)(\omega)= f(\omega), \quad \omega \in S.$
\item The range of $T_p$ is dense in $\mathcal L^p(\partial \Omega).$
\item
$$
W^{(1)}_{p,0}(\Omega)= \left\{ f(x)\in W^{(1)}_p(\Omega): (T_pf)(\omega)=0 \right\}.
$$
\end{enumerate} 
Properties (1), (2) and (3) are the main three properties that make useful the standard trace maps in the case of Lipschitz boundaries. This shows that our trace maps effectively replace the standard trace maps in the case where the boundaries are not Lipschitz.  

We consider now a simple application of our trace theorem to the Dirichlet problem. Following Section 2 of Chapter 8 of \cite{gt} we formulate the Dirichlet problem as follows. Let $\Omega$ be any bounded domain in $\mathbf R^n.$ Let $L$ be the formal differential operator,
\beq\label{2.17}
Lf:=D_i \left( a^{ij}(x) D_j f+ b^i(x) f\right)+  c^i(x) D_i f+ d(x)f,
\ene
with  $D_i:= \frac{\partial}{\partial x_i}, i=1,\cdots,n,$ and where we use the convention of summation over repeated indices. We assume that the coefficients $a^{ij}, b^i, c^i$ and $d, i,j=1,\cdots n$ are Lebesgue measurable and bounded in $\Omega,$ that $ a^{ij}$ is strongly elliptic, i.e., that for some $ \lambda >0,$
\beq\label{2.18}
a^{ij}(x)\, \zeta_i\, \zeta_j \geq \lambda |\zeta|^2, \quad \forall x \in \Omega, \zeta \in \mathbf R^n.
\ene
Furthermore, we assume the following negativity condition
\beq\label{2.19}
\int_{\Omega}\, \left( d  \varphi- b^i D_i \varphi \right)\, dx \leq 0, \quad \forall \varphi \in C^\infty_0(\Omega), \varphi \geq 0.
\ene
Moreover, let $ g, h^i, i=1,\cdots n,$ be locally integrable functions in $\Omega.$ Formally the Dirichlet problem can be written as,
\beq\label{2.19.b}
\begin{array}{l}
Lf(x)= g(x)+ D_i h^i(x), \quad x \in \Omega,\\
f(x)= h(x), \quad  x \in \partial \Omega, 
\end{array}
\ene
for some function $h$ defined in $\partial \Omega.$ In a precise mathematical sense the Dirichlet problem is formulated as follows. By a weak (or generalized) solution to the equation,

\beq\label{2.20}
Lf(x)= g(x)+ D_i h^i(x), \quad x \in \Omega,
\ene
we mean a function $ f \in W^{(1)}_2(\Omega)$ such that, 
$$
\int_{\Omega}\,\left[ \left( a^{ij} D_j f+ b^i f\right)\, D_i \varphi - (c^i D_i f + d \,f )\, \varphi\right]\, dx=\int_{\Omega}\, \left(h^i D_i \varphi-g \varphi \right)\, dx \quad     \forall \varphi \in C^1_0(\Omega).
$$

 We have the following theorem

\begin{theorem}\label{dp}
Suppose that the functions $a^{ij}, b^i, c^i, 1 \leq i,j \leq n$ and $d$ are Lebesgue measurable and bounded in $\Omega.$ Further, assume that \eqref{2.18} and \eqref{2.19} hold, and  that the functions $ g, h^i \in L^2(\Omega), i=1,\cdots,n.$ Then, for every $h\in W^{(1)}_2(\Omega) $, the equation \eqref{2.20} has a unique weak solution $ f \in W^{(1)}_2(\Omega)$ such that, $f-h \in  W^{(1)}_{2,0}(\Omega).$
\end{theorem} 

\noindent {\it Proof: }This result is Theorem 8.3 in page 181 of \cite{gt} where the proof is given.

\qed

Theorem \ref{dp} is a general existence and uniqueness result.  However, it gives no information in how the solution $f$ takes the values of 
$h$ in $\partial \Omega.$ In fact, the situation is even worse, because as we only know that  $f, h \in W^{(1)}_2(\Omega),$ in the absence of a trace map, we can not even give a mathematical meaning to the values   of $f$ and $h$ in $\partial \Omega.$ In the next theorem we provide a precise answer to this issue using our trace map. 
\begin{theorem}
Assume  that the domain $\Omega$ is star shaped with respect to  the origin with a continuos function that characterizes the boundary. 
Further, suppose that the functions $a^{ij}, b^i, c^i, 1 \leq i,j \leq  n$ and $d$ are Lebesgue measurable and bounded in $\Omega.$ Moreover, suppose that \eqref{2.18} and \eqref{2.19} hold, and  that the functions $ g, h^i \in L^2(\Omega), i=1,\cdots,n.$ Then, for every $h\in W^{(1)}_2(\Omega) $, the unique weak solution $ f \in W^{(1)}_2(\Omega)$ to \eqref{2.20}  such that, $f-h \in  W^{(1)}_{2,0}(\Omega),$ satisfies 
$T_2f=T_2 h$ in $\mathcal L^2(\partial \Omega).$ In particular,  $f(\omega)= h(\omega)$ for $m_{\partial \Omega}$- almost every $\omega \in \partial \Omega.$ Equivalently,   $f( b(\nu)\nu)= h(b(\nu)\nu)$ for Lebesgue almost every $\nu  \in S_{n-1}.$ 
\end{theorem}
\noindent{\it Proof:} The theorem follows from Theorem \ref{dp} and \eqref{2.16.b} since as $ f-h \in  W^{1}_{2,0}(\Omega),$ we have that
$T_2f= T_2 h.$
\qed

Being able to give a precise mathematical meaning to the values of the data and the solution on the boundary, and proving that they coincide up to a set of measure zero, is an important piece of information, not only from the mathematical point of view, but also from the numerical  perspective and furthermore,  in order to use the solution to our problem in concrete applications.

 \section{The case of $\mathbf R^n$}\sss
We consider now traces on  star-shaped surfaces   of functions in $W^{(s)}_2(\mathbf R^n),n\geq 2 $. We already discussed this problem in \cite{we}. However, there we studied the traces within the context of the spectral and scattering theory of strongly propagative systems of equations and in the particular case of the slowness surface of these systems. Here we present our results in all its generality  and in a more detailed way.  

 We define  star shaped {\it hypersurfaces} as in the case of the boundary of a bounded domain, but we only assume that  the function $b(\nu)$ is Lebesgue measurable in $S_{n-1}.$
\begin{definition} \label{def3.1}
{\rm We say that a set   $S$ in $\mathbf R^n$   is a star shaped with respect to the origin  {\it hypersurface}  characterized by a Lebesgue measurable function,   if there is a function $b(\nu) >0 $ defined for $ \nu \in S_{n-1}$ that is Lebesgue measurable in $S_{n-1},$ such that,}
\beq\label{3.1} 
S=\left\{ x \in \mathbf R^n: |x|= b\left(\frac{x}{|x|}\right)\right\}.
\ene
Furthermore, we suppose that there is a constant $ C_1>0$ such that, $ 0 <C_1 <b(\nu) < 1/C_{1},$ for  $ \nu \in S_{n-1}.$  
\end{definition}
Let us define the Fourier transform on $L^2(\mathbf R^n)$ as follows,
$$
(F_nf)(k):= \frac{1}{(2\pi)^{n/2}}\, \int_{\mathbf R^n}\, e^{-ik\cdot x}\, f(x)\, dx.
$$
Recall that, for $l=1,\cdots,$
\beq\label{3.1.a}
k_j^l \, (F_nf)(k)= \left( F_n\, (-i)^l  \, \frac{\partial^l}{\partial x^l_j}  f(x)\right)(k), \quad j=1,\cdots, n,
\ene
and that,
\beq\label{3.1.b}
 \frac{\partial^l}{\partial k_j^l} (F_nf)(k)=  \left( F_n(-i)^l\, x_j^l  f(x)\right)(k), \quad j=1,\cdots, n.
\ene

As usual,  the Sobolev space $ W^{(s)}_2(\mathbf R^n), s >0,$ is defined as the set of all functions $ f(x)\in L^2(\mathbf R^n)$ such that the  Fourier transform $(F_nf)(k)$ satisfies,
$(1+k^2)^{s/2}\, (F_nf)(k)\in L^2(\mathbf R^n),$ with norm given by,
\beq\label{3.1.c}
\left\| f \right\|_{W^{(s)}_2(\mathbf R^n)}:=  \left\| (1+k^2)^{s/2}\, F_n f\right\|_{L^2(\mathbf R^n)}.
\ene
Clearly, $ \mathcal S$ is dense on $W^{(s)}_2(\mathbf R^n).$ Observe that when $s=1$  the norm \eqref{3.1.c} is equivalent to the norm \eqref{2.3} with $ \Omega= \mathbf R^n$ and $p=2.$ We use the norm \eqref{3.1.c} for the simplicity of notation.

We define a measure in $S$ and the space $\mathcal L^2(S)$  as in the case of star shaped domains in Section \ref{sec2}.
We first  define the measure in $S$. Let $\mathcal P_S$ be the following  function from $S$ onto $S_{n-1}.$ 
\beq\label{3.2}
\mathcal P_S(\omega):= \nu= \frac{\omega}{|\omega|}, \,\text{for}\, \omega \in S.
\ene 
The inverse function is given by,
\beq\label{3.3}
\mathcal P^{-1}_S(\nu):= \omega= b(\nu)\, \nu, \, \text{for}\, \nu \in S_{n-1}.
\ene
As before $\mathcal P_S$ is onto and one-to-one.  
\begin{definition}\label{def3.2}{\rm
We say that a set $ O\subset S, $ is $ S$-measurable if $ \mathcal P_S(O)$ is Lebesgue measurable in $S_{n-1}.$ Further, for any  
$S$-measurable set $O,$ we define its measure,  $ m_{S}(O),$ by
\beq\label{3.4} 
m_{S}(O):= m_{S_{n-1}}(\mathcal P_S(O)).
\ene
\label{3.5}
 By $\mathcal M_{S}$ we denote the $\sigma$-algebra of all $S$-measurable sets. Note that  $m_{S}$ is a $\sigma$-additive measure on   $\mathcal M_{S}.$}
\end{definition}
As before, for a function,  $f(\omega),$ defined on $S,$ we designate  by $f_{S_{n-1}}(\nu)$ the function,
\beq
f_{S_{n-1}}(\nu):= f(b(\nu)\,  \nu), \, \nu \in S_{n-1}.
\ene
Then,  for any set $ A \subset \mathbf R,$
$$
\mathcal P_S\left(f^{-1}(A)\right)= f^{-1}_{S_{n-1}}( A).
$$
We have  that  $f(\omega)$ is $m_{S}$-measurable if and only if $f_{S_{n-1}}(\nu)$ is Lebesgue measurable.
For any $ m_{S}$-integrable function $f(\omega),$ we have that,
\beq\label{3.6}
\int_{S}\, f(\omega)\, dm_{S}(\omega)= \int_{S_{n-1}}\,f_{S_{n-1}}(\nu)\, dm_{S_{n-1}}(\nu).
\ene

\begin{definition} \label{def3.3} {\rm
We denote by $\mathcal L^2(S) $ the Hilbert space of all complex valued, $m_{S}$-measurable functions, $ f(\omega)$ such that
$|f(\omega)|^2$ is integrable over $S$ with respect to the measure $m_{S}.$The scalar product in   $\mathcal L^2(S)$  is given by
\beq\label{3.7}
\ds \int_{S}\, f(\omega) \, \overline{g(\omega)}\, dm_{S}(\omega)=\ds \int_{S_{n-1}}\,f_{S_{n-1}}(\nu)\, \ds  \overline{ g_{S_{n-1}}(\nu)}\,dm_{S_{n-1}}(\nu).
\ene
}
\end{definition}

 We find it convenient to define the following weighted space.
\begin{definition}\label{def3.4}
{\rm We denote by $\mathcal L^2_b(S)$ the Hilbert space of all complex valued $m_S$-measurable functions that are square integrable on $S$ with the scalar product,
$$
\ds (f,g)_{\mathcal L^2_b(S)}:= \ds\int_{S}\, f(\omega)\, \overline{g(\omega)}\, b^n(\omega/|\omega|)\, dm_{S}(\omega)= \int_{S_{n-1}}\,f_{S_{n-1}}(\nu)\, 
      \overline{g_{S_{n-1}} (\nu)}\,\,   b^n(\nu)\,dm_{S_{n-1}}(\nu).
      $$ }
\end{definition}
As there is a constant $C_1$ such that, $ 0 < C_1 < |b(\nu)| <1/C_1, \nu \in S_{n-1},$ the norms of $\mathcal L^2(S)$ and of $\mathcal L^2_b(S)$ are equivalent.

The basic idea to prove our trace theorem is to parametrize  $x \in\mathbf R^n$ as $ x= \rho\, \omega, \rho > 0, \omega \in S.$ In this parametrization taking the restriction to $S$ means to take a sharp value of $\rho,$ and this should not require regularity in $\omega.$ Hence, assuming that $b(\nu)$ is Lebesgue measurable should be sufficient. As we will show this is actually true.

We prepare the following results that we need for the proof of our trace theorem. 
Suppose  that  $ f(x) \in C^\infty_0(\mathbf R^n).$ Then, using spherical coordinates, we have that,
\beq\label{3.8}
\int_{\mathbf R^n}\, |f(x)|^2 \, dx= \int_{S_{n-1}}\, dm_{S_{n-1}}(\nu)\, \int_0^\infty\,dr\, r^{n-1}\,  |f(r\nu)|^2 \,=
 \int_{S_{n-1}}\, dm_{S_{n-1}}(\nu)\,  h(\nu),
 \ene
 where
 \beq\label{3.8.a}
 h(\nu):= \int_0^\infty\, dr\, r^{n-1}\, |f(r\nu)|^2.
 \ene
 For each fixed $ \nu \in S_{n-1},$ we perform the change of variable  $ \rho= r/ b(\nu)$ in the one dimensional integral on the right-hand side of \eqref{3.8.a} and we obtain that,
 \beq\label{3.8.b}
 h(\nu)=   b^n(\nu)\ \int_0^\infty\, d\rho\, \rho^{n-1}\,  |f(\rho\, b(\nu)\,\nu)|^2.
 \ene 
Further, by \eqref{3.8}, \eqref{3.8.b},
\beq\label{3.8.c}
\int_{\mathbf R^n}\, |f(x)|^2 \, dx= \int_{S_{n-1}}\, dm_{S_{n-1}}(\nu)\, b^n(\nu)\, \int_0^\infty\, d\rho\, \rho^{n-1}\, |f(\rho\, b(\nu)\,\nu)|^2.
\ene
We define the function,
\beq\label{3.8.d}
g(\rho,\nu):= \rho^{n-1}\, |f(\rho\, b(\nu)\,\nu)|^2.
\ene
By \eqref{3.8.c} and \eqref{3.8.d}, it follows that,
\beq\label{3.8.e}
\int_{\mathbf R^n}\, |f(x)|^2 \, dx= \int_{S_{n-1}}\, dm_{S_{n-1}}(\nu)\, b^n(\nu)\, \int_0^\infty\, d\rho\, g(\rho,\nu).
\ene
 Since  $ dm_{S_{n-1}}(\nu) \times d\rho$ is a product measure, by   Fubini's theorem we can exchange the order of the integration of the function $g(\rho,\nu)$ in the right-hand side of
 \eqref{3.8.e}, to obtain that,
 \beq\label{3.8.f}
  \int_0^\infty\, d\rho\ \int_{S_{n-1}}\, dm_{S_{n-1}}(\nu)\, b^n(\nu)\, g(\rho,\nu)= \int_{S_{n-1}}\, dm_{S_{n-1}}(\nu)\, b^n(\nu)\, \int_0^\infty\, d\rho\, g(\rho,\nu).
\ene
By \eqref{3.8.e} and \eqref{3.8.f} 
 \beq\label{3.8.g}
 \int_{\mathbf R^n}\, |f(x)|^2 \, dx= \int_0^\infty\, d\rho\ \int_{S_{n-1}}\, dm_{S_{n-1}}(\nu)\, b^n(\nu)\, g(\rho,\nu).
 \ene
 
 Let $U$ be the following operator, first defined for $ f(x)\in C^\infty_0(\mathbf R^n),$
 \beq\label{3.10}
 (Uf)(\rho, \omega):= \rho^{\frac{n-1}{2} }f(\rho \,\omega), \quad \rho >0, \omega\in S.
 \ene
 Then, by \eqref{3.6}
 \beq\label{3.10.b}
\ds \int_0^\infty\, d\rho\, \int_{S}\, | (Uf)(\rho,\omega)|^2\, b^n(\omega/|\omega|) dm_{S}(\omega)= \int_0^\infty\, d\rho\,\rho^{n-1}\, \int_{S_{n-1}}\, b^n(\nu)\,   |f(\rho\, b(\nu) \,\nu )|^2\, dm_{S_{n-1}}(\nu).
\ene
By \eqref{3.8.d}, \eqref{3.8.g} and \eqref{3.10.b}
\beq\label{3.10.c}
\ds \int_0^\infty\, d\rho\, \int_{S}\, | (Uf)(\rho,\omega)|^2\, b^n(\omega/|\omega|) dm_{S}(\omega)= \int_{\mathbf R^n}\, |f(x)|^2 \, dx.
\ene

 Hence, $U$ extends to a unitary operator from $L^2(\mathbf R^n)$ onto $ L^2((0,\infty);\mathcal L^2_b(S)).$ The unitarity of $U$ plays an important role in the spectral and scattering theory of strongly propagative systems of equations in \cite{we}.
 
  let $h(\lambda)$ be a function in $C^\infty(\mathbf R)$ such that,
 $ 0 \leq h (\lambda)\leq 1, h(\lambda)= 0,  \lambda < 1 /4, h(\lambda)=1, 1/2 < \lambda < \infty$. For $\rho >0$ let us denote $h_\rho(\lambda)= h(\lambda/\rho).$   Note that, $ h_\rho(\lambda)=0,  \lambda < \rho/4, h_\rho(\lambda)=1, \rho/2 <  \lambda < \infty.$ For $ f(x)\in \mathcal S$ we define the operator,
 \beq\label{3.9.b}  
 \left(L_\rho(z)f\right)(x,\omega)= (1+x^2)^{z/2}\, \frac{1}{\sqrt{2\pi}}\,\int_{\mathbf R}\, d \lambda\, e^{-i x\, \lambda}\, h_\rho(\lambda)\, \left(F_n^{-1} (1+k^2)^{-z/2}\, F_nf\right)(\lambda \,\omega),
 \ene
 where $ x \in \mathbf R, \omega \in S, z= \alpha+i\beta \in \mathbf C, 0 \leq \alpha \leq 4.$ 
 
 For any pair of Banach spaces, $X,Y,$ we denote by $\mathcal B(X,Y)$ the Banach space of all bounded operators from $X$, into $Y.$
 \begin{lemma} \label{lemm3.1}
 The operator $ L_\rho(z), 0 \leq \text{Re}\, z \leq 4,$ extends to a bounded operator from $L^2(\mathbf R^n)$ into $L^2(\mathbf R; \mathcal L^2(S)).$ Furthermore, for any $\rho_0 >0,$ and any $0 \leq \alpha \leq 4,$ there is a constant, $C_{\rho_0}(\alpha),$ that depends only on $\rho_0,$ and $\alpha,$ such that,
 
 \beq\label{3.9.c}
 \left\|  L_{\rho}(\alpha+i\beta)\right\|_{\mathcal B \left(L^2(\mathbf R^n), L^2(\mathbf R; \mathcal L^2(S))\right) }  \leq C_{\rho_0}(\alpha)\, \frac{1}{\rho^{(n-1)/2}}, \quad \rho \geq \rho_0, 0 \leq \alpha \leq 4.
 \ene 
 \end{lemma}
 
 \noindent{\it Proof:}
 Denote,
 $$
 \psi(\lambda,\omega):=  \left(F_n^{-1} (1+k^2)^{-i\beta/2}\, F_nf\right)(\lambda \,\omega).
$$
 Using Parseval's identity for the one dimensional Fourier transform, taking into account that for some constant $C_1,$ $ 0< C_1  < b(\omega/|\omega|) < 1/C_{1}, \omega \in S,$  as in the support of $h_\rho(\lambda),  \lambda > \rho/4,$  and using \eqref{3.10.b}, we have that,
 \beq\label{3.10.c.1}\begin{array}{l}
\int_{S}\, dm_S(\omega)\, \int_{\mathbf R}\, dx    \left| (L_\rho(i \beta) f)(x,\omega)\right|^2=\int_{S}\, dm_S(\omega)\, \int_{\rho/4}^\infty\, d \lambda\, h_{\rho}(\lambda)^2\, |\psi(\lambda \,\omega)|^2 \leq \\\\
C \, \frac{1}{\rho^{n-1}}\,\int_{S}\, dm_S(\omega)\,  b^n(\omega/|\omega|)\,  \int_0^\infty \, d\lambda\, \lambda^{n-1} |\psi(\lambda \,\omega)|^2= C\, \frac{1}{\rho^{n-1}}\, \left\| UF_n^{-1} (1+k^2)^{-i\beta/2}\, F_nf\right\|^2_{ L^2((0,\infty); \mathcal L^2_b(S))}=  \\\\
C\,  \frac{1}{\rho^{n-1}}\,\left\| f\right\|^2_{L^2(\mathbf R^n)},
\end{array}
\ene
where in the last equality we used  that $U$ is unitary and Parseval's identity for the Fourier transform on $ \mathbf R^n.$
It follows that $T_\rho(i\beta), \beta \in \mathbf R,$ extends to a bounded operator from $L^2(\mathbf R^n)$ into $L^2(\mathbf R; \mathcal L^2(S))$ and, 
\beq\label{3.11.0}
  \left\|L_\rho(i \beta)f\right\|_{L^2(\mathbf R; \mathcal L^2(S))} \leq C \, \frac{1}{\rho^{(n-1)/2}}\, \left\|  f\right\|_{L^2(\mathbf R^n)}.
  \ene
Furthermore, for $f$ in Schwartz space denote,
\beq\label{3.11}
(M_\rho(z)f)(x,\omega):= (1+x^2)^{-z/2}\, (L_\rho(z)f)(x, \omega). 
\ene
Then,
\beq\label{3.12}
\left\|L_\rho(4+i \beta)f\right\|_{L^2(\mathbf R; \mathcal L^2(S))} \leq C \left\|M_\rho(4+i \beta)f\right\|_{L^2(\mathbf R; \mathcal L^2(S))} + C \left\|x^4 M_\rho(4+i \beta)f\right\|_{L^2(\mathbf R; \mathcal L^2(S))}. 
\ene
We denote,
\beq\label{3.14}
\varphi( \lambda, \omega):= \left(F_n^{-1} (1+k^2)^{-(4+i\beta)/2}\, F_nf\right)(\lambda \,\omega).
\ene

Arguing as \eqref{3.10.c.1}, we prove that,
\beq\label{3.13}\begin{array}{l}
 \left\|M_\rho(4+i \beta)f\right\|^2_{L^2(\mathbf R; \mathcal L^2(S))}  \leq  C \, \frac{1}{\rho^{n-1}}\,\int_{S}\, dm_S(\omega)\,  b^n(\omega/|\omega|)\,  \int_0^\infty \, d\lambda\, \lambda^{n-1} |\varphi(\lambda \,\omega)|^2 \leq  C\,\frac{1}{\rho^{(n-1)}}\, \left\| f\right\|^2_{L^2(\mathbf R^n)}.
\end{array}
\ene

Moreover, using \eqref{3.1.a} with $n=1,$ with  $k$ replaced by $x$ and $x$ replaced by $\lambda,$ we obtain that,
\beq\label{3.15}
\begin{array}{l}
\left\|x^4 M_\rho(4+i \beta)f\right\|^2_{L^2(\mathbf R; \mathcal L^2(S))}  \leq C\ \, \frac{1}{\rho^{n-1}}\,  \int_{S}\, dm_S(\omega)\,  b^n(\omega/|\omega|)\,  \int_0^\infty \, d\lambda\, \lambda^{n-1}  
  \left|\frac{\partial^4}{\partial \lambda^4}\,  h_\rho(\lambda) \varphi(\lambda,\omega)\right|^2.
 \end{array} \ene
Furthermore,
\beq\label{3.16}
\frac{\partial}{\partial \lambda } \varphi(\lambda,\omega) = i  \sum_{l=1}^n \omega_j\, \frac{1}{(2\pi)^{n/2}}\,\int_{\mathbf R^n}\,dk \,e^{i \lambda \omega \cdot k} \, k_j (1+k^2)^{-(4+i\beta)/2}\,  \,  (F_n f)(k).
\ene 
Then, using  \eqref{3.10.b},  \eqref{3.16} as well as its derivatives with respect to $\lambda,$ and Parseval's identity on $L^2(\mathbf R^n),$ we obtain that,
\beq\label{3.17}
  \int_{S}\, dm_S(\omega)\,  b^n(\omega/|\omega|)\,  \int_0^\infty \, d\lambda\, \lambda^{n-1} 
  \left|\frac{\partial^4}{\partial \lambda^4}\,  h_\rho(\lambda) \varphi(\lambda \,\omega)\right|^2  \leq C  \,   \left\|  f\right\|^2_{L^2(\mathbf R^n)}.
  \ene
  Further, by \eqref{3.12}, \eqref{3.13}, \eqref{3.15} and \eqref{3.17}
  \beq\label{3.18}
  \left\|L_\rho(4+i \beta)f\right\|_{L^2(\mathbf R; \mathcal L^2(S))} \leq C  \, \frac{1}{\rho^{(n-1)/2}}\,  \left\|  f\right\|_{L^2(\mathbf R^n)}.
  \ene
 Note that it follows from our estimates that the constants $C$ in \eqref{3.11.0} and \eqref{3.18} are uniform   for $\rho  \geq \rho_0,$ for any $\rho_0 >0.$   By \eqref{3.11.0}, \eqref{3.18} and Hadamard three lines  theorem,   \cite{rs2}, page 33,  $L_\rho(\alpha+i \beta), 0 \leq \alpha \leq 4, $ is  bounded   from $L^2(\mathbf R^n)$ into $L^2(\mathbf R ;\mathcal L^2(S)),$ and \eqref{3.9.c} holds.
 \qed
\begin{remark}\label{rem3.1}
Suppose that $ s > 1/2.$ Then, for any $\rho,\rho_1>0 $ there is a constant $C_s$ that depends only on $s$ such that,
\beq\label{3.18.0}
\int_{\mathbf R}\, dx\, \frac{|e^{ix \rho}- e^{ix\rho_1}|^2}{(1+ |x|)^{2s}}\leq C_s \left\{\begin{array}{l} |\rho-\rho_1|^{2s-1}, \quad \frac{1}{2} < s < \frac{3}{2},\\\\
 |\rho-\rho_1|^2\, (1+ |\ln|\rho-\rho_1||),\,\quad s=\frac{3}{2}, \\\\
|\rho-\rho_1|^2, \quad s > \frac{3}{2}.
\end{array}\right.
\ene
\end{remark}
\noindent{\it Proof:}
We  can assume that $|\rho-\rho_1|\leq1$. Then, For any $R >1,$ 
$$\begin{array}{l}
\ds\int_{\mathbf R}\, dx\, \frac{|e^{ix \rho}- e^{ix\rho_1}|^2}{(1+ |x|)^{2s}}=  \int_{|x|\geq R-1 }\, dx\, \frac{|e^{ix \rho}- e^{ix\rho_1}|^2}{(1+ |x|)^{2s}}+
\int_{|x|\leq R-1}\, dx\, \frac{|e^{ix \rho}- e^{ix\rho_1}|^2}{(1+ |x|)^{2s}}\leq\\\\ C_s
\ds\left\{\begin{array}{l}
 R^{1-2s}+ R^{3-2s}\, |\rho-\rho_1|^2, \quad  \frac{1}{2} < s < \frac{3}{2}, \\\\
\ds R^{1-2s}+ |\rho-\rho_1|^2\,  \ln R,\quad s= 3/2,\\\\
R^{1-2s}+|\rho-\rho_1|^2, \quad s > \frac{3}{2}.
\end{array}\right. 
\end{array}
$$ 
Finally, \eqref{3.18.0} follows taking $ R=|\rho-\rho_1|^{-1}.$
\qed

In the following theorem we state our trace map.
\begin{theorem}\label{theo3.1}
Let  the set  $S$ in $\mathbf R^n$ be  a star shaped with respect to the origin  {\it hypersurface} characterized by a Lebesgue  measurable function. Then, for every $ \rho >0$  and every $ s >1/2,$ there is a trace map $T_s(\rho)$ that is bounded from 
$ W^{(s)}_2(\mathbf R^n)$ into $\mathcal L^2(S)$ such that for every function $f(x)$ in the space of Schwartz,
\beq\label{3.18.b}
(T_s(\rho)f)(\omega)= f(\rho \,\omega).
\ene
Moreover, for every $\rho_0>0$ there is constant $C_{\rho_0}$ such that,
\beq\label{3.18.c}
\left\| T_s(\rho) \right\|_{\mathcal B\left(W^{(s)}_2(\mathbf R^n), \mathcal L^2(S)\right)} \leq C_{\rho_0}\, \frac{1}{\rho^{(n-1)/2}}, \quad 
\rho \geq \rho_0,
\ene
and, furthermore,
\beq\label{3.18.d}\begin{array}{l}
\left\| T_s(\rho)-T_s(\rho_1) \right\|_{\mathcal B\left(W^{(s)}_2(\mathbf R^n), \mathcal L^2(S)\right)} \leq C_{\rho_0}\, \ds\frac{1}{\min[\rho^{(n-1)/2},\rho_1^{(n-1)/2}]} \times \\\\
   \left\{\begin{array}{l} |\rho-\rho_1|^{s-1/2}, \quad \frac{1}{2} < s < \frac{3}{2},\\\\
 |\rho-\rho_1|\,  (1+ \sqrt{|\ln|\rho-\rho_1||}),\,\quad s=\frac{3}{2}, \\\\
|\rho-\rho_1|, \quad s > \frac{3}{2},
\end{array}\right.
\end{array}
\ene
for $ \rho,\rho_1 \geq \rho_0.$
Moreover, the range of $T(\rho)$ is dense in $\mathcal L^2(S), $
\beq\label{3.18.bb}
\overline{T(\rho) W^{(s)}_2(\mathbf R^n)}= \mathcal L^2(S).
\ene
\end{theorem}

\noindent{\it Proof:}   For  $ f(x)\in  \mathcal S$ we define,
\beq\label{3.18.e}
 (T(\rho)\, f)(\omega):= f(\rho \,\omega), \quad \omega \in S.
\ene
Let us denote,
\beq\label{3.18.f} 
g(x):= \left(F_n^{-1} (1+k^2)^{s/2}\, (F_nf)(k)\right)(x).
\ene
Note that,
\beq\label{3.18.g} 
f(\rho\, \omega)= \frac{1}{\sqrt{2\pi}}\, \int_{\mathbf R}\, dx \,e^{i\rho \,x}\,  (1+x^2)^{-s/2} \left(L_\rho(s) g\right)(x,\omega).
\ene
 Then, by Schwarz inequality, for $ s> 1/2,$
 $$
 |f(\rho\, \omega)|^2 \leq C \int_{\mathbf R}\, dx\,  \left|(L_\rho(s) g)(x,\omega)\right|^2.
 $$
 Hence, for any $ \rho_0 >0,$ and for any $ \rho \geq \rho_0,$
  \beq\label{3.19}\begin{array}{l}
 \int_{S}\,dm_{S}(\omega)  \, |f(\rho\, \omega)|^2  \leq
  C  \int_{S}\,dm_{S}(\omega)\, \,\int_{\mathbf R}\, dx\,  \left|(L_\rho(s) g)(x,\omega)\right|^2  = \\\\
  C\, \|L_\rho(s)\,g\|^2_{L^2(\mathbf R; \mathcal L^2(S))}
\leq C \,C_{\rho_0}(s) \, \frac{1}{\rho^{n-1}}\  \|g\|_{L^2(\mathbf R^n)}^2= C\, C_{\rho_0}(s)\frac{1}{\rho^{n-1}}\, \|f\|^2_{W^{(s)}_2(\mathbf R^n)},
\end{array}
\ene
where we used  \eqref{3.9.c}. It follows that $T(\rho)$ is bounded and hence, it extends uniquely to a bounded operator from $W^{(s)}_2(\mathbf R^n)$ into $\mathcal L^2(S)$ and  \eqref{3.18.c} holds. Let us now prove \eqref{3.18.d}. Without loss of generality we can assume that $ \rho_1 > \rho \geq \rho_0.$ Suppose that $f(x)\in \mathcal S.$  Since $h_\rho(\lambda)= 1,
\lambda > \rho/2,$ we have that $ h_\rho(\rho_1)=1.$ In consequence, it follows that,
\beq\label{3.20} 
 (T(\rho_1)f)(\omega)=f(\rho_1\, \omega)= \frac{1}{\sqrt{2\pi}}\, \int_{\mathbf R}\, dx\, e^{i\rho_1 x}\,  (1+x^2)^{-s/2} \left(L_\rho(s) g\right)(x,\omega),
\ene
where $g(x)$ is defined in \eqref{3.18.f}. Further, by \eqref{3.18.e}, \eqref{3.18.g} and \eqref{3.20}
\beq\label{3.21}
(T(\rho)f)(\omega)- (T(\rho_1)f)(\omega)= f(\rho\,\omega )-f(\rho_1\, \omega)= \frac{1}{\sqrt{2\pi}}\, \int_{\mathbf R}\, dx \left(e^{i\rho \,x}- e^{i \rho_1x}\right)\,  (1+x^2)^{-s/2} 
\left(L_\rho(s) g\right)(x,\omega).
\ene
Hence, by Schwarz inequality, 
\beq\label{3.22}\begin{array}{l}
\ds\left\| (T(\rho)-T(\rho_1)) f \right\|^2_{\mathcal L^2( S)}  \leq \frac{1}{\sqrt{2\pi}} \,\int_{ S}\, dm_S(\omega)\, \int_{\mathbf R}\, dx\, \frac{|e^{ix \rho}- e^{ix\rho_1}|^2}{(1+ |x|)^{2s}}  \int_{\mathbf R}\, dx\, \left|\left(L_\rho(s) g\right)(x,\omega)\right|^2 = \\\\
\ds   \frac{1}{\sqrt{2\pi}} \, \int_{\mathbf R}\, dx\, \frac{|e^{ix \rho}- e^{ix\rho_1}|^2}{(1+ |x|)^{2s}}\, \|L_\rho(s)g\|^2_{L^2(\mathbf R;\mathcal L^2(S) )} \leq   \frac{1}{\sqrt{2\pi}} \,\int_{\mathbf R}\, dx\, \frac{|e^{ix \rho}- e^{ix\rho_1}|^2}{(1+ |x|)^{2s}} C^2_{\rho_0}(s) \frac{1}{\rho^{n-1}}\,  \left\|f \right\|^2_{W^{(s)}_2}(\mathbf R^n),       
\end{array}
\ene
where we used \eqref{3.9.c} and \eqref{3.18.f}. Then, \eqref{3.18.d} follows from  \eqref{3.18.0} and \eqref{3.22}.

We proceed now to prove \eqref{3.18.bb} as in the case of Theorem \ref{theo2.4}. Take any $f(\omega) \in \mathcal L^2(S).$ It follows that, $f_{S}(\nu):= f(b(\nu) \nu) \in L^2(S_{n-1}).$ As $C^\infty(S_{n-1})$ is dense on $L^2(S_{n-1})$ there is sequence  $g_m(\nu) \in C^\infty(S_{n-1})$ such that,
\beq\label{3.23}
\lim_{m\to \infty}\| f_{S_{n-1}}-g_m \|_{L^2(S_{n-1})}=0.
\ene
Let $C_1$ be such that, that $ 0 < C_1 <b(\nu) < 1/C_{1}, \nu \in S_{n-1}.$ Take $ h(r) \in C^\infty_0((0,\infty)),$ such that $ h(r)=0, 0< r <   C_1/4, h(r)=1, C_1/ 2 < r < 2/ C_{1}, h(r)=0, 
r > 3 /C_{1},$ and denote, $h_\rho(r):= h(r/\rho).$ We define,
\beq\label{3.24}
f_m(x):= h_\rho(|x|)\, g_m\left(\frac{x}{|x|}\right) \in  C^\infty_0(\mathbf R^n) \subset W^{(s)}_2(\mathbf R^n).
\ene
Finally,
\beq\label{3.25}  \begin{array}{l}
\lim_{m \to \infty}\|f -T(\rho)  f_m  \|^2_{\mathcal L^2(S)}=\lim_{m \to \infty} \int_{S_{n-1}}\, |f(b(\nu)\nu)- h_\rho( \rho\, b(\nu))\, g_m(\nu)|^2\, dm_{S_{n-1}}(\nu)= \\\\
 \lim_{m \to \infty} \int_{S_{n-1}}\, |f(b(\nu)\nu)- g_m(\nu)|^2\, dm_{S_{n-1}}(\nu)= \lim_{m \to \infty} \int_{S_{n-1}}\, |f_{S_{n-1}}(\nu)- g_m(\nu)|^2\, dm_{S_{n-1}}(\nu)=0.
\end{array} \ene
\qed

In Theorem  A.1 of \cite{we} we proved a result like Theorem ~\ref{theo3.1} in the case where $S$ is the slowness surface of a strongly propagative system of  equations. A slowness  surface is defined in terms of a continuous function, $\lambda(k), k \in \mathbf R^n\setminus\{0\}$  with values in $(0,\infty)$ that is homogeneous of order one, i.e. $\lambda(\rho\, k)= \rho\, \lambda(k),  \rho >0, k \in \mathbf R^n\setminus \{0\}.$ The slowness surface is given by $S:=\{ k \in \mathbf R^n: \lambda(k)=1\}.$
This corresponds in Definition \ref{def3.1}  to $ b(\nu)= 1/\lambda(\nu), \nu \in S_{n-1}.$ Note that in Theorem A.1 of \cite{we}  the precise estimates \eqref{3.18.c}, \eqref{3.18.d} and the density in $\mathcal L^2(S)$ of the range of the trace map $T(\rho)$  are not proven. In \cite{we} the trace map in Theorem A.1 is applied to the spectral and the scattering theory of the strongly propagative systems of equations. Note that the H\"older continuity of the trace map plays an essential role in the spectral and  scattering theory of strongly propagative systems of equations.  Furthermore, the boundedness of the operator $T_s(\rho)$  was proven in Proposition 1.11 in page 115 of \cite{ya}. In fact, estimate  (1.35) in page 115 of \cite{ya}, with  $\kappa=1,$ is equivalent to the boundedness of  the operator $T_s(\rho)$ that we prove in Theorem ~\ref{theo3.1}. 
Remark that in Proposition 1.11 in page 115 of \cite{ya}  a notation that is different from ours is used. Moreover, in  Proposition 1.11 in page 115 of \cite{ya} the H\"older continuity and the precise estimates  \eqref{3.18.c} and \eqref{3.18.d} for the norm of $T_s(\rho)$ are not obtained. The proof of Proposition 1.11 in page 115 of \cite{ya}  is different from ours. It is based in the Mourre method. For the use of the Mourre method in this context see \cite{wemou}.

\section{The coarea formula}\sss
For the  coarea formula   in the case where the level sets are given by a Lipschitz function see \cite{fe}. For the case where  the level sets are given by a function in a Sobolev space look to   \cite{msz}. See also  \cite{fr}.    Let us consider the coarea formula in a simple case.  Denote,
$$
a(x):= \frac{|x|} {b(x/|x|)},    
$$
and suppose that $b(\nu), \nu \in S_{n-1}$ is smooth and that for some $ C_1 >0, C_1 < b(\nu) < 1/C_{1}, \nu \in S_{n-1}$. Denote,
$$
S_\lambda=\{ x \in \mathbf R^n : a(x)=\lambda \},
$$
and
$$
\Omega_\lambda:= \{ x \in \mathbf R^n : a(x) < \lambda\}.
$$
Note that $S$ in \eqref{3.1} satisfies $S=S_1.$
Denote by $V_\lambda$ the volume of $\Omega_\lambda$ with respect to the Lebesgue measure. Then, by the coarea formula,
\beq\label{4.1}
V_\lambda= \int_0^\lambda\, d \rho\, \int_{S_\rho}\, \frac{1}{| \nabla a(x)|}\, d\beta_\rho(x),
\ene
where $d\beta_\rho(x)$ denotes the Lebesgue measure of  $S_\rho.$
Furthermore,
\beq\label{4.2}
\frac{d}{d \lambda} V(\lambda)=  \int_{S_\lambda}\, \frac{1}{| \nabla a(x)|}\, d\beta_\lambda ( x).
\ene
The right hand-side of \eqref{4.2} is a trace with the standard trace map in the surface $S_\lambda.$ 
 
  We now prove that with our trace map we can generalize \eqref{4.1} and \eqref{4.2} to the case where $b(\nu)$ is only assumed to be Lebesgue measurable and, of course to satisfy $ 0 < C < b(\nu) < 1/C, \nu \in S_{n-1}.$
Using spherical coordinates, we have that,
\beq\label{4.3}\begin{array}{l}
V_\lambda = \int_{S_{n-1}}\, dm_{S_{n-1}}(\nu)\, \int_0^{ \lambda b(\nu)}\,dr\, r^{n-1}=\\\\
\ds\frac{\lambda^n}{n}\,   \int_{S_{n-1}}\, dm_{S_{n-1}}(\nu)\, b^n(\nu)=
\int_0^\lambda\, d \rho \, \rho^{n-1}\,\int_{S_{n-1}}\, \, b^n(\nu)dm_{S_{n-1}}(\nu).
\end{array}\ene
 As in the case of $S=S_1$ we  define a measure on $S_\rho$  as follows.  By $\mathcal P_{S_\rho}$  we denote  the function from $S_\rho$ onto $S_{n-1}$ given by, 
$$
\mathcal P_{S_\rho}(\omega):= \nu= \frac{\omega}{|\omega|}, \,\text{for}\, \omega \in S_{\rho}.
$$ 

A set $ O\subset S_{\rho}, $ is $ S_\rho$-measurable if $ \mathcal P_{S_\rho}(O)$ is Lebesgue  measurable in $S_{n-1},$ and the measure of $O$ is given by,
$$
m_{S_\rho}(O):=  \rho^{n-1}\,m_{S_{n-1}}(\mathcal P_{S_\rho}(O)).
$$
As in Section 2, for a function,  $f(\omega),$ defined on $S_\rho,$ we denote,
$$
f_{S_{n-1}}(\nu):= f( \rho \,b(\nu)\,  \nu), \, \nu \in S_{n-1},
$$
and   $f(\omega)$ is $S_\rho$-measurable if and only if $f_{S_{n-1}}(\nu)$ is Lebesgue measurable.
For any $ S_\rho$-integrable functions $f(\omega),$ we have that,
\beq\label{3.6.b}
\int_{S_\rho}\, f(\omega)\, dm_{S_\rho}(\omega)=  \rho^{n-1}\, \int_{S_{n-1}}\,f_{S_{n-1}}(\nu)\, dm_{S_{n-1}}(\nu).
\ene

 With this definition \eqref{4.3} reads,
 \beq\label{4.4.a}
 V_\lambda= \int_0^\lambda\, d \rho \, \int_{S_{\rho}}\,   b^n(\omega/|\omega|)\,  dm_{S_{\rho}}(\omega).
 \ene
 Furthermore, by \eqref{4.4.a}
 \beq\label{4.4.b}
 \frac{d}{d\lambda} V_\lambda= \int_{S_{\lambda}}\,    b^n(\omega/|\omega|)\, dm_{S_{\lambda}}(\omega).
 \ene
Formulae  \eqref{4.4.a} and \eqref{4.4.b} generalize, respectively, \eqref{4.1} and \eqref{4.2} to the case where $b(\nu)$ is only assumed to be Lebesgue measurable.  

Moreover, in the case where $b(\nu)$ is smooth, for any function $f(x)\in C^\infty_0(\mathbf R^n),$ the coarea formula implies that,
 \beq\label{4.5}
   \int_{\mathbf R^n}\, f(x) \, dx= \int_0^\infty\, d \rho\, \int_{S_\rho}\,\,f(x)\, \frac{1}{| \nabla a(x)|}\, d\beta_\rho(x) .
   \ene
   With our method, we prove as in the proof of \eqref{3.8.g} that,
   \beq\label{4.6}
   \int_{\mathbf R^n}\, f(x)\, dx= \int_0^\infty\, d\rho\, \rho^{n-1} \int_{S_{n-1}}\, \, f(\rho\, b(\nu)\,\nu)\,\, b^n(\nu)\,
dm_{S_{n-1}}(\nu).
   \ene
 Then, using \eqref{3.6.b}  and \eqref{4.6} we obtain that
    \beq\label{4.7}
   \int_{\mathbf R^n}\, f(x)\, dx= \int_0^\infty\, d\rho\,  \int_{S_\rho}\,\, f(\omega)\,  b^n(\omega/|\omega|)\, dm_{S_{\rho}}(\omega).
   \ene 
   Formula \eqref{4.7} generalizes the coarea formula \eqref{4.5} to the case where $b(\nu)$ is only assumed to be Lebesgue measurable. Note that with the notation of Theorem \ref{theo3.1} with $S=S_1$, equation \eqref{4.7} can be written as follows,
   \beq\label{4.7.a.a}
    \int_{\mathbf R^n}\, f(x)\, dx= \int_0^\infty\, d\rho\,     \rho^{n-1}\, \int_{S}\, (T(\rho) f)(\omega)\,  b^n(\omega/|\omega|)\, dm_{S}(\omega).
   \ene 
   Since $T(\rho)$ is bounded, \eqref{4.7.a.a} extends by continuity to functions in $L^1(\mathbf R^n)\cap W^{(s)}_2(\mathbf R^n), s >1/2.$
   
   Actually, it can be verified directly that our formulae coincide with the ones given by the coarea formula when $b(\nu)$ is smooth. For simplicity, we do the calculations for $n=2,$  but the result is also true for $ n\geq 3$ with a similar computation. 
As in Section \ref{sec2}   let us take take polar coordinates $ x_1= r \cos\theta, x_2= r \sin\theta, $ with $ r=|x|, 0 \leq \theta < 2 \pi,$ and we denote $ b(\theta):= b(\nu(\theta)),$   $ e(\theta):=1/ b(\theta),$ with $\nu(\theta)$ as in \eqref{unitvec}.   In these coordinates,
$$
 a(r,\theta):= a(r \nu(\theta))= r \,e(\theta).
$$
 Then, on $S_{\rho},$
\begin{equation}\label{4.7.b}
|\nabla a(r,\theta)|= \sqrt{(e(\theta))^2+(e'(\theta))^2}.
\end{equation}
Let us now write $S_{\rho}$ in parametric form,
$$
S_{\rho}= \{ x \in \mathbf R^2: x= \varphi(\theta), 0 \leq \theta \leq 2 \pi\}.
$$
Recall that, $e_{x_1},e_{x_2}$ are, respectively, the unit vectors along the $x_1$ and the $x_2$ axis. Then, in polar coordinates,
$$
\varphi(\theta)= \frac{\rho}{e(\theta)}\,(\cos\theta e_{x_1}+ \sin\theta e_{x_2}).
$$
It follows that,
$$
\left| \varphi'(\theta) \right|= \frac{\rho}{(e(\theta))^2}\, \sqrt{(e(\theta))^2+ (e'(\theta))^2}.
$$
Then,
$$
d\beta_\rho(x)\equiv d\beta_\rho(\theta)= \frac{\rho}{e^2(\theta)}\, \sqrt{(e(\theta))^2+ (e'(\theta))^2}\, d\theta.
$$
Finally, using \eqref{4.7.b}  and that $e(\theta) = 1/ b(\theta),$ we obtain that, 
\beq\label{4.8}
\ds \frac{1}{|\nabla a(x)|}\, d\beta_\rho(x) = \rho\, b^2(\theta)\, d\theta.
\ene
By \eqref{4.8} the right-hand sides of \eqref{4.1} and \eqref{4.4.a} are the same, the right-hand sides of \eqref{4.2} and \eqref{4.4.b} coincide, and  the right-hand sides of   \eqref{4.5} and \eqref{4.7} take the same value.

\noindent{\bf Acknowledgement}This paper was partially written while I was visiting 
INRIA Saclay \^Ile-de-France and ENSTA. I thank Anne-Sophie Bonnet-BenDhia  and Patrick Joly for their kind hospitality.  I thank Vladimir Maz'ya for his detailed information on the literature on trace maps and on the coarea formula.

\end{document}